\documentclass{elsart}
\usepackage{amsfonts}
\usepackage{amsmath}
\usepackage{amssymb}
\usepackage{graphicx}
\journal{Topology}
\begin{document}

%TCIMACRO{\TeXButton{Begin frontmatter}{\begin{frontmatter}}}
%BeginExpansion
\begin{frontmatter}
%EndExpansion
%

%TCIMACRO{\TeXButton{Title}{\title
%{A combinatorial Yamabe flow in three dimensions}}}%
%BeginExpansion
\title{A combinatorial Yamabe flow in three dimensions}%
%EndExpansion
%

%TCIMACRO{\TeXButton{Author}{\author{David Glickenstein}}}
%BeginExpansion
\author{David Glickenstein}%
%EndExpansion
%

%TCIMACRO{\TeXButton{Address}{\address
%{Department of Mathematics, University of Arizona, Tucson, AZ 85721, USA}}}%
%BeginExpansion
\address
{Department of Mathematics, University of Arizona, Tucson, AZ 85721, USA}%
%EndExpansion
%

%TCIMACRO{\TeXButton{Begin abstract}{\begin{abstract}}}%
%BeginExpansion
\begin{abstract}
%EndExpansion
A combinatorial version of Yamabe flow is presented based on Euclidean
triangulations coming from sphere packings. The evolution of curvature is then
derived and shown to satisfy a heat equation. The Laplacian in the heat
equation is shown to be a geometric analogue of the Laplacian of Riemannian
geometry, although the maximum principle need not hold. It is then shown that
if the flow is nonsingular, the flow converges to a constant curvature metric.
%TCIMACRO{\TeXButton{End abstract}{\end{abstract}}}%
%BeginExpansion
\end{abstract}
%EndExpansion
%

%TCIMACRO{\TeXButton{Begin keywords}{\begin{keyword}}}%
%BeginExpansion
\begin{keyword}
%EndExpansion
curvature flow, Yamabe flow, sphere packing, Laplacian, discrete Riemannian
geometry
%TCIMACRO{\TeXButton{End keyword}{\end{keyword}}}
%BeginExpansion
\end{keyword}
%EndExpansion
%

%TCIMACRO{\TeXButton{End frontmatter}{\end{frontmatter}}}%
%BeginExpansion
\end{frontmatter}
%EndExpansion

\section{Introduction}

In his proof of Andreev's theorem in \cite{thurston13}, Thurston introduced a
conformal geometric structure on two-dimensional simplicial complexes which is
an analogue of a Riemannian metric. He then used a version of curvature to
prove the existence of circle packings (see also Marden-Rodin
\cite{mardenrodin} for more exposition). Techniques very similar to elliptic
partial differential equation techniques were used by Y. Colin de Verdi\`{e}re
\cite{verdiere} to study conformal structures and circle packings. Cooper and
Rivin in \cite{cooperrivin} then defined a version of scalar curvature on
three-dimensional simplicial complexes and used it to look at rigidity of
sphere packings along the lines of Colin de Verdi\`{e}re.

Inspired by this work, Chow and Luo \cite{chowluo} defined several
combinatorial Ricci flows on two-dimensional simplicial complexes, one for
each constant curvature model space. They were able to show that the flows
converge to constant curvature if a circle packing exists whose nerve is the
one-skeleton of the triangulation. The reader is also directed to some later
work of Luo on how these flows evolve the conformal structure \cite{luoyamabe}
. We shall use Cooper and Rivin's combinatorial scalar curvature to define
combinatorial Yamabe flow on three-dimensional simplicial complexes which is a
three-dimensional analogue of Chow and Luo's work when the triangles are
modeled on Euclidean triangles. We shall look at the evolution of curvature
from a geometric viewpoint, understanding the heat equation on curvature which
is induced by the flow. The flow turns out not to be parabolic in the usual
sense of Laplacians on graphs, an analytic property which we study in a
related paper \cite{glickensteinmaxprinciple}. The geometric flow perspective
is very much inspired by Richard Hamilton's works on the Ricci flow, e.g.
\cite{hamilton:singularities}.

The combinatorial Yamabe flow is a way of studying prescribed scalar curvature
on simplicial complexes, which we might call the combinatorial Yamabe problem.
The Yamabe problem has been studied in great detail (see \cite{leeparker} for
a good overview). The Riemannian case has been solved by the work of Aubin
\cite{aubin} and Schoen \cite{schoen}. Yamabe flow in the smooth category has
been studied by Hamilton and others. We refer the reader to R. Hamilton
\cite{hamilton:surfaces} and R. Ye \cite{ye}.

\section{Geometric structures and combinatorial manifolds}

\begin{figure}
[ptb]
\begin{center}
\includegraphics[
natheight=3.555200in,
natwidth=3.555200in,
height=1.7269in,
width=1.7269in
]
{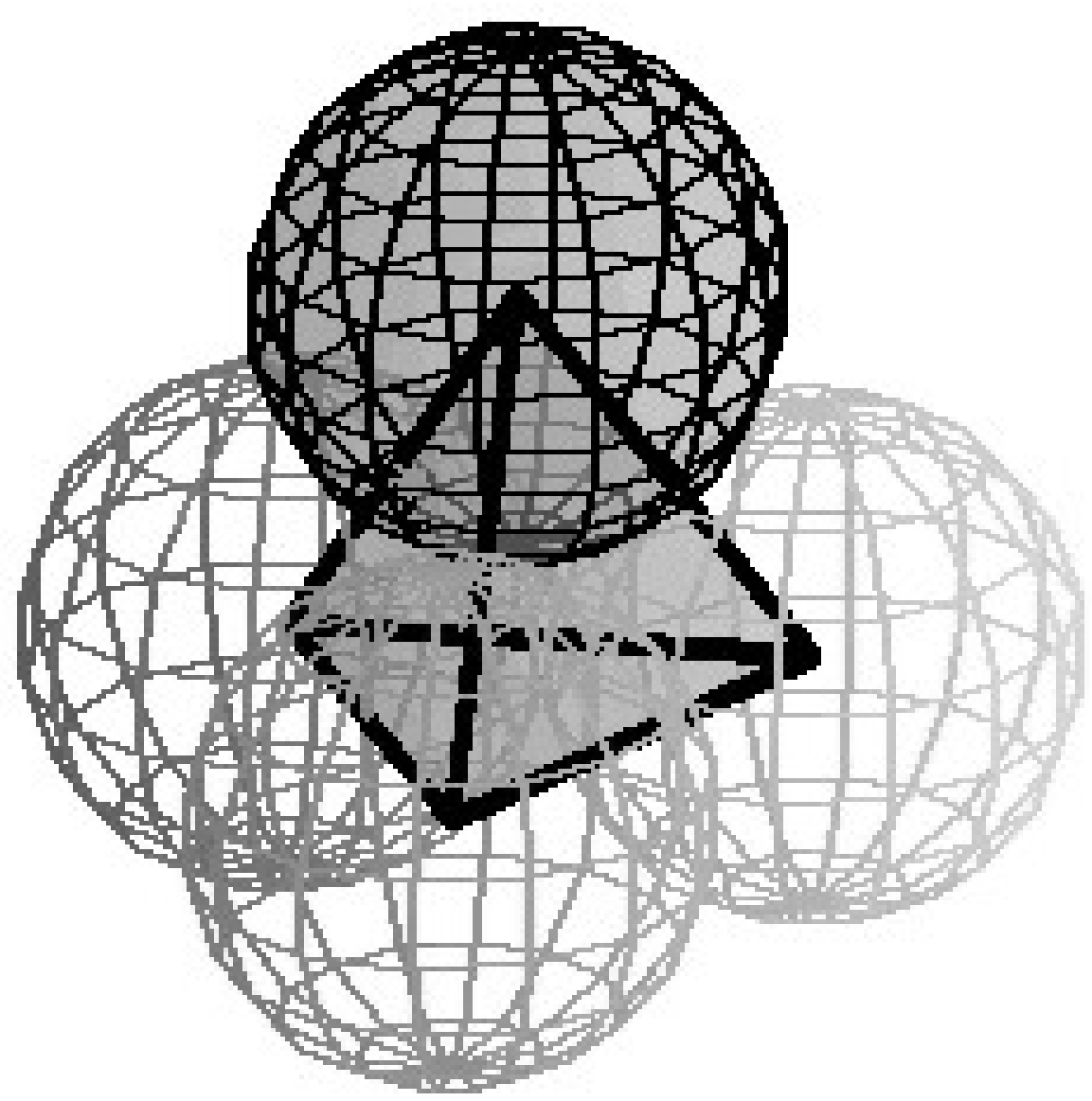}
\caption{Tetrahedron with balls at the vertices.}
\label{tetrawithballs}
\end{center}
\end{figure}

We essentially take our formalism from Cooper-Rivin in \cite{cooperrivin}. We
shall use the notation $f_{i}$ to denote evaluation of a function $f$ at $i$
in a finite set and $f\left(  t\right)  $ to denote evaluation at $t$ in an
interval. Let $\mathcal{S=}\left\{  \mathcal{S}_{0},\mathcal{S}_{1}
,\ldots,\mathcal{S}_{n}\right\}  $ be a simplicial complex of dimension $n,$
where $\mathcal{S}_{i}$ is the $i$-dimensional skeleton. We define a metric
structure as a map
\[
r:\mathcal{S}_{0}\mathcal{\rightarrow}\left(  0,\infty\right)
\]
such that for every edge $\left\{  i,j\right\}  \in\mathcal{S}_{1}$ between
vertices $i$ and $j,$ the length of the edge is $\ell_{ij}=r_{i}+r_{j}.$ Any
such metric structure, or a particular tetrahedron within such a structure, is
called conformal and the set of all is called the conformal class. We can
think of this as having an $n$-dimensional sphere packing whose nerve is the
collection of edges $\mathcal{S}_{1},$ although it is not necessarily an
actual sphere packing. One such conformal tetrahedron is shown in Figure
\ref{tetrawithballs}.

\begin{figure}
[ptbptb]
\begin{center}
\includegraphics[
natheight=5.926800in,
natwidth=5.926800in,
height=1.8422in,
width=1.8422in
]
{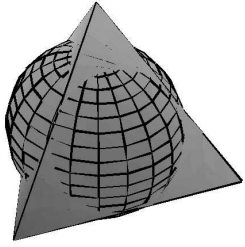}
\caption{Tetrahedron with circumscripted sphere.}
\label{tetra with circumscripted sphere}
\end{center}
\end{figure}

Conformal tetrahedra are also called circumscriptible tetrahedra, and the
condition on the edges is equivalent to the condition that there exists a
sphere tangent to each of the edges of the tetrahedron \cite[Chapter
9.B.1]{altshiller} (we call this sphere the \emph{circumscripted sphere} since
it is circumscripted by the tetrahedron) as seen in Figure
\ref{tetra with circumscripted sphere}.

The function $r$ determines the 2-dimensional faces since there is a
one-to-one correspondence between triples $\left(  r_{i},r_{j},r_{k}\right)  $
and triples of sides for Euclidean triangles given by
\[
r_{i}=\frac{1}{2}\left(  \ell_{ij}+\ell_{ik}-\ell_{jk}\right)
\]
and so forth. We shall also put the restriction that each higher dimensional
simplex can be realized as a Euclidean simplex. We shall return to this
condition a little later.

Each metric of this type is in some sense conformal to the metric where all
$r_{i}=1,$ since they can be gotten by rescaling the function $r$ at each
point. This is similar to multiplying a Riemannian metric $g$ by a function
$f^{2}$ at every point to get a new metric $f^{2}g$ which is conformal to the
metric $g.$ The metric structure $\left\{  r_{i}\right\}  _{i\in
\mathcal{S}_{0}}$ determines the geometry, which comes from the lengths of the
edges, similar to the way a Riemannian metric determines the metric space
structure of a Riemannian manifold.

\begin{figure}
[ptb]
\begin{center}
\includegraphics[
natheight=4.301500in,
natwidth=6.844700in,
height=1.82in,
width=2.8892in
]
{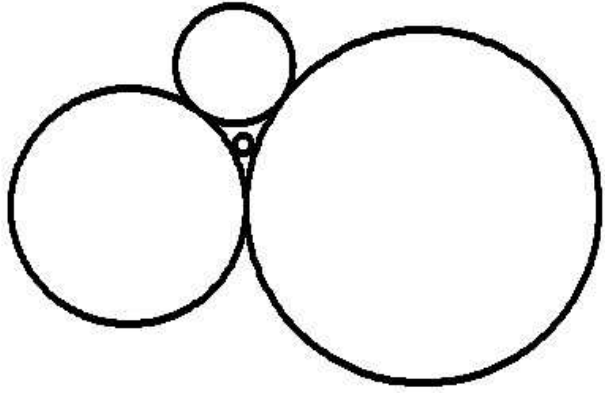}
\caption{Failure of four circles to be mutually tangent.}
\label{circlediagram}
\end{center}
\end{figure}

In the sequel we shall limit ourselves primarily to three dimensions. Cooper
and Rivin \cite[Section 3]{cooperrivin} observe that for a collection
$\left\{  r_{i},r_{j},r_{k},r_{\ell}\right\}  $ to determine a Euclidean
tetrahedron, we can use Descartes' circle theorem, also called Soddy's
theorem, which says that four circles in the plane of radii $r_{i},r_{j}
,r_{k},r_{\ell}$ are externally tangent if
\[
Q_{ijk\ell}\doteqdot\left(  \frac{1}{r_{i}}+\frac{1}{r_{j}}+\frac{1}{r_{k}
}+\frac{1}{r_{\ell}}\right)  ^{2}-2\left(  \frac{1}{r_{i}^{2}}+\frac{1}
{r_{j}^{2}}+\frac{1}{r_{k}^{2}}+\frac{1}{r_{\ell}^{2}}\right)  =0.
\]
For a nice proof of Soddy's theorem, see \cite{pedoe}. We also direct the
reader to the interesting article \cite{lagariasmallowswilks} on Descartes'
circle theorem. Looking at the proof it is clear that if this quantity is
negative, then we get three circles which are mutually tangent and a circle in
the middle which cannot be tangent to all the others, as seen in Figure
\ref{circlediagram},
and hence we cannot form a Euclidean tetrahedron from spheres of these radii.
If $Q_{ijk\ell}$ is positive, we can form a Euclidean tetrahedron
corresponding to $\left\{  r_{i},r_{j},r_{k},r_{\ell}\right\}  .$ So our
condition for nondegeneracy of the tetrahedron is
\[
Q_{ijk\ell}=\left(  \frac{1}{r_{i}}+\frac{1}{r_{j}}+\frac{1}{r_{k}}+\frac
{1}{r_{\ell}}\right)  ^{2}-2\left(  \frac{1}{r_{i}^{2}}+\frac{1}{r_{j}^{2}
}+\frac{1}{r_{k}^{2}}+\frac{1}{r_{\ell}^{2}}\right)  >0.
\]
We call $Q_{ijk\ell}$ the nondegeneracy quadratic. As noted in \cite[Section
793]{altshiller}, $Q_{ijk\ell}$ is really $4$ divided by the square of the
radius of the sphere tangent to each of the edges (the circumscripted sphere),
and is related to the volume in the following way:
\[
V_{ijk\ell}^{2}=\frac{1}{9}r_{i}^{2}r_{j}^{2}r_{k}^{2}r_{\ell}^{2}Q_{ijk\ell}.
\]
Thus if we consider the formula for the square of the volume as formal, the
nondegeneracy condition is that $V_{ijk\ell}^{2}>0.$

Now we shall define a quantity $K$ called the curvature. For a Euclidean
tetrahedron with vertices $\left\{  i,j,k,\ell\right\}  $ we define the solid
angle $\alpha_{ijk\ell}$ at a vertex $i$ as the area of the triangle on the
unit sphere cut out by the planes determined by $\left\{  i,j,k\right\}
,\left\{  i,j,\ell\right\}  ,\left\{  i,k,\ell\right\}  $ where $i$ is the
center of the sphere. Note that the solid angle is also sometimes called the
trihedral angle. If we define $\beta_{ijk\ell}$ as the dihedral angle in the
tetrahedron $\left\{  i,j,k,\ell\right\}  $ along the edge $\left\{
i,j\right\}  ,$ which is also an angle of the aforementioned spherical
triangle, the Gauss-Bonnet theorem gives the formula for the solid angle as
\[
\alpha_{ijk\ell}=\beta_{ijk\ell}+\beta_{ikj\ell}+\beta_{i\ell jk}-\pi.
\]
Note that solid angles $\alpha_{ijk\ell}$ are symmetric in the last three
indices and dihedral angles $\beta_{ijk\ell}$ are symmetric in the first two
and in the last two indices. We can now define the curvature $K_{i}$ at a
vertex $i$ as
\[
K_{i}\doteqdot4\pi-\sum_{\left\{  i,j,k,\ell\right\}  \in\mathcal{S}_{3}
}\alpha_{ijk\ell}.
\]
Note that the sum is over $j,k,\ell$ since the vertex $i$ is fixed.

The curvature $K_{i}$ can be thought of as a scalar curvature since it
measures the difference at a given vertex between the total angles in
Euclidean space and the total angles of the complex. It was initially looked
at by Cooper-Rivin in \cite{cooperrivin}. They found that a metric structure
cannot be deformed (staying conformal) while keeping the scalar curvature constant.

Constant scalar curvature is a critical point of the total curvature
functional
\[
T=\sum K_{i}r_{i}.
\]
Cooper-Rivin showed that the space of nondegenerate simplices is not convex,
but that the function $T$ is weakly concave as a function of the $r_{i}$ and
strongly concave if the condition $\sum r_{i}=1$ is imposed. We cannot use
this to show that there is a unique constant scalar curvature metric on a
given conformal class, but any constant scalar curvature metric is a local
minimum of $T.$

Another way to prove this is to look at the functional $T$ as a function of
$s_{i}=1/r_{i}.$ Then the set of nondegenerate simplices is convex since it is
the intersection of a cone with a half-space. However, upon computing the
Hessian of $T$ we can only prove that the function is concave in a
neighborhood of constant curvature. Much like scalar curvature in the smooth
case, there may be several constant scalar curvature metrics in a given
conformal class. However, we have yet to find a complex admitting two constant
scalar curvature metrics.

\begin{rem}
\label{vertex transitive remark}Just because a topological manifold admits a
constant sectional curvature $0$ metric, a given triangulation of that
manifold may not admit a metric with curvature $0.$ For any vertex transitive
triangulation, i.e. a triangulation such that the same number of tetrahedra
meet at every vertex, the triangulation with $r_{i}=1$ for all $i\in$
$\mathcal{S}_{0}$ is constant curvature. If $d$ is the degree of the vertex,
i.e. the number of tetrahedra meeting at each vertex, the curvature must be
$4\pi-d\left(  3\cos^{-1}\left(  1/3\right)  -\pi\right)  .$ Hence in order
for the curvature to be zero, $d=\left(  4\pi\right)  /\left(  3\cos
^{-1}\left(  1/3\right)  -\pi\right)  \approx\allowbreak22.\,\allowbreak795.$
Thus no vertex transitive triangulation of the torus admits a zero curvature
metric. This observation is similar to the one noted in
\cite{amborncarforamarzuoli} about dynamical triangulations.
\end{rem}

\begin{rem}
The curvature considered here is different than the one considered in the
Regge calculus by Regge \cite{regge} and others (for instance, \cite{hamber},
\cite{hamberwilliams}, \cite{frohlich}). Our curvature is concentrated at the
vertices while Regge's curvature is concentrated at the edges in dimension 3.
The solid angle which we use is reminiscent of the interpretation of Ricci
curvature as solid angle deficit (see, for instance, the introduction of
\cite{besse}).
\end{rem}

\section{Combinatorial Yamabe flow}

We now define combinatorial Yamabe flow on the metric structure as
\begin{equation}
\frac{dr_{i}}{dt}=-K_{i}r_{i}\label{combinatorial yamabe flow}
\end{equation}
for each $i\in\mathcal{S}_{0}.$ Note how similar this looks to the Yamabe flow
on Riemannian manifolds, which is
\[
\frac{\partial}{\partial t}g_{ij}=-Rg_{ij}
\]
where $g_{ij}$ is the Riemannian metric and $R$ is its scalar curvature. In
particular, both preserve their respective conformal class. We use the term
`combinatorial' since this is used by Chow and Luo, but it is really more of a
piecewise linear or piecewise Euclidean flow because it depends on the
geometry of the triangulation and not just the topological, or combinatorial,
structure. Much like the Yamabe flow and Ricci flow, the evolution of
curvature will play a key role in understanding the behavior of this equation.
Next we shall compute this evolution.

Recall the Schl\"{a}fli formula, which, for a Euclidean tetrahedron denoted by
the complex $\left\{  \mathcal{T}_{0},\mathcal{T}_{1},\mathcal{T}
_{2},\mathcal{T}_{3}\right\}  ,$ gives that
\[
\sum_{\left\{  i,j\right\}  \in\mathcal{T}_{1}}\ell_{ij}d\beta_{ijk\ell}=0
\]
(see Milnor \cite{milnorschlafli} for a proof). We can reorganize this as
Cooper-Rivin \cite{cooperrivin} do to get
\[
\sum_{i\in\mathcal{T}_{0}}r_{i}\left(  d\beta_{ijk\ell}+d\beta_{ikj\ell
}+d\beta_{i\ell jk}\right)  =\sum_{i\in\mathcal{T}_{0}}r_{i}\,d\alpha
_{ijk\ell}=0
\]
or
\[
r_{i}\frac{\partial\alpha_{ijk\ell}}{\partial r_{i}}+r_{j}\frac{\partial
\alpha_{jik\ell}}{\partial r_{i}}+r_{k}\frac{\partial\alpha_{kij\ell}
}{\partial r_{i}}+r_{\ell}\frac{\partial\alpha_{\ell ijk}}{\partial r_{i}}=0.
\]
Since there are only four vertices in $\mathcal{T}_{0},$ we may denote the
solid angle at vertex $i$ by $\alpha_{i}$ without fear of confusion. Then we
consider
\begin{equation}
A\doteqdot\sum_{i\in\mathcal{T}_{0}}r_{i}\alpha_{i} \label{define A}
\end{equation}
so
\begin{align*}
dA  &  =\sum_{i\in\mathcal{T}_{0}}\alpha_{i}\,dr_{i}+\sum_{i\in\mathcal{T}
_{0}}r_{i}\,d\alpha_{i}\\
&  =\sum_{i\in\mathcal{T}_{0}}\alpha_{i}\,dr_{i}
\end{align*}
by the Schl\"{a}fli formula. We thus have
\[
\frac{\partial A}{\partial r_{i}}=\alpha_{i}
\]
and hence
\[
\frac{\partial\alpha_{i}}{\partial r_{j}}=\frac{\partial^{2}A}{\partial
r_{i}\partial r_{j}}=\frac{\partial\alpha_{j}}{\partial r_{i}}
\]
by commuting the partial derivatives. In our expanded notation, which we shall
use for complexes larger than one simplex, this says
\[
\frac{\partial\alpha_{ijk\ell}}{\partial r_{j}}=\frac{\partial\alpha_{jik\ell
}}{\partial r_{i}}.
\]
Using our derivation from the Schl\"{a}fli formula we also have
\begin{equation}
r_{i}\frac{\partial\alpha_{ijk\ell}}{\partial r_{i}}+r_{j}\frac{\partial
\alpha_{ijk\ell}}{\partial r_{j}}+r_{k}\frac{\partial\alpha_{ijk\ell}
}{\partial r_{k}}+r_{\ell}\frac{\partial\alpha_{ijk\ell}}{\partial r_{\ell}
}=0. \label{schlafli}
\end{equation}
This equality has a much more geometric interpretation; it says that the
directional derivative of the angle $\alpha_{ijk\ell}$ in the direction of
scaling $\left(  r_{i},r_{j},r_{k},r_{\ell}\right)  $ is zero, since if these
are scaled equally, the new tetrahedron is similar to the original and hence
all angles remain the same.

The evolution of curvature is
\begin{align*}
\frac{d}{dt}K_{i}  &  =-\sum_{\left\{  i,j,k,\ell\right\}  \in\mathcal{S}_{3}
}\frac{d}{dt}\alpha_{ijk\ell}\\
&  =-\sum_{\left\{  i,j,k,\ell\right\}  \in\mathcal{S}_{3}}\left(
\frac{\partial\alpha_{ijk\ell}}{\partial r_{i}}\frac{dr_{i}}{dt}
+\frac{\partial\alpha_{ijk\ell}}{\partial r_{j}}\frac{dr_{j}}{dt}
+\frac{\partial\alpha_{ijk\ell}}{\partial r_{k}}\frac{dr_{k}}{dt}
+\frac{\partial\alpha_{ijk\ell}}{\partial r_{\ell}}\frac{dr_{\ell}}{dt}\right)
\\
&  =\sum_{\left\{  i,j,k,\ell\right\}  \in\mathcal{S}_{3}}\left(
\frac{\partial\alpha_{ijk\ell}}{\partial r_{i}}K_{i}r_{i}+\frac{\partial
\alpha_{ijk\ell}}{\partial r_{j}}K_{j}r_{j}+\frac{\partial\alpha_{ijk\ell}
}{\partial r_{k}}K_{k}r_{k}+\frac{\partial\alpha_{ijk\ell}}{\partial r_{\ell}
}K_{\ell}r_{\ell}\right) \\
&  =\sum_{\left\{  i,j,k,\ell\right\}  \in\mathcal{S}_{3}}\left(
\frac{\partial\alpha_{ijk\ell}}{\partial r_{j}}r_{j}\left(  K_{j}
-K_{i}\right)  +\frac{\partial\alpha_{ijk\ell}}{\partial r_{k}}r_{k}\left(
K_{k}-K_{i}\right)  \right. \\
&  \quad\quad\quad\quad\quad+\left.  \frac{\partial\alpha_{ijk\ell}}{\partial
r_{\ell}}r_{\ell}\left(  K_{\ell}-K_{i}\right)  \right)
\end{align*}
using (\ref{schlafli}). We call the coefficients
\[
\Omega_{ijk\ell}\doteqdot\frac{\partial\alpha_{ijk\ell}}{\partial r_{j}}
r_{j}.
\]
In this notation we see that the evolution of curvature is
\[
\frac{d}{dt}K_{i}=\sum_{\left\{  i,j,k,\ell\right\}  \in\mathcal{S}_{3}
}\left[  \Omega_{ijk\ell}\left(  K_{j}-K_{i}\right)  +\Omega_{ikj\ell}\left(
K_{k}-K_{i}\right)  +\Omega_{i\ell jk}\left(  K_{\ell}-K_{i}\right)  \right]
.
\]

In order to compute the coefficients $\Omega_{ijk\ell}$ we need to compute the
partial derivatives of the solid angles. We do this computation using the
following formulas from Euclidean geometry. Recall that $\alpha_{ijk\ell}$
refers to the solid angle of tetrahedron $\left\{  i,j,k,\ell\right\}  $ at
$i$ and that $\beta_{ijk\ell}$ refers to the dihedral angle of tetrahedron
$\left\{  i,j,k,\ell\right\}  $ along the edge $\left\{  i,j\right\}  .$ We
also need the face angles. Denote the angle of the triangle $\left\{
i,j,k\right\}  $ at the vertex $i$ by $\gamma_{ijk}.$ We can then use the law
of cosines and the expression for area in terms of sines to compute
$\gamma_{ijk}$ as
\begin{align}
\cos\gamma_{ijk}  &  =\frac{\ell_{ij}^{2}+\ell_{ik}^{2}-\ell_{jk}^{2}}
{2\ell_{ij}\ell_{ik}}=\frac{r_{i}^{2}+r_{i}r_{j}+r_{i}r_{k}-r_{j}r_{k}
}{\left(  r_{i}+r_{j}\right)  \left(  r_{i}+r_{k}\right)  }
\label{cosinefaceangle}\\
\sin\gamma_{ijk}  &  =\frac{2A_{ijk}}{\ell_{ij}\ell_{ik}}=\frac{2\sqrt
{r_{i}r_{j}r_{k}\left(  r_{i}+r_{j}+r_{k}\right)  }}{\left(  r_{i}
+r_{j}\right)  \left(  r_{i}+r_{k}\right)  }\nonumber
\end{align}
where $\ell_{ij}=r_{i}+r_{j}$ is the length of edge $\left\{  i,j\right\}  $
and $A_{ijk}=\sqrt{r_{i}r_{j}r_{k}\left(  r_{i}+r_{j}+r_{\ell}\right)  }$ is
the area of triangle $\left\{  i,j,k\right\}  $ by Heron's formula.

Using the law of cosines for the face angles (\ref{cosinefaceangle}), we can
compute the evolution of the face angles, which turns out to be
\[
\frac{d}{dt}\gamma_{ijk}=-2\frac{A_{ijk}}{P_{ijk}}\left(  \frac{K_{j}-K_{i}
}{\ell_{ij}}+\frac{K_{k}-K_{i}}{\ell_{ik}}\right)
\]
where we have introduced the notation $P_{ijk}=2\left(  r_{i}+r_{j}
+r_{k}\right)  $ for the perimeter of the triangle $\left\{  i,j,k\right\}  .$
It should be noted that this computation was entirely formal, and is thus the
same formula derived by Chow and Luo \cite{chowluo} for simplicial surfaces.
If we define the curvature of a surface to be $k_{i}=2\pi-\sum\gamma_{ijk}$
then the formula for evolution of the face angles implies that the evolution
of curvature is in fact parabolic in the usual sense of Laplacians on graphs
(this is studied in greater detail in \cite{glickensteinmaxprinciple}). The
curvature evolution turns out to be
\[
\frac{dk_{i}}{dt}=\left(  \Delta k\right)  _{i}
\]
where the Laplacian is the one defined by Z. He in \cite{he}. We shall explore
this aspect in the next section.

The face angles are used to compute the dihedral angles and solid angles via
spherical geometry. If we consider the solid angle formed by three planes, say
those determined by $\left\{  i,j,k\right\}  ,\left\{  i,j,\ell\right\}  ,$
and $\left\{  i,k,\ell\right\}  ,$ we see that the planes intersect the sphere
and form a spherical triangle. It is clear that the angles of this triangle
are the dihedral angles $\beta_{ijk\ell},\beta_{ikj\ell},\beta_{i\ell jk}$ and
that the length of the sides of this triangle are the face angles
$\gamma_{ijk},\gamma_{ij\ell},\gamma_{ik\ell},$ hence the relationship between
the dihedral angles and the face angles can be expressed in terms of the
spherical law of cosines, which says
\begin{equation}
\cos\beta_{ijk\ell}=\frac{\cos\gamma_{ik\ell}-\cos\gamma_{ijk}\cos
\gamma_{ij\ell}}{\sin\gamma_{ijk}\sin\gamma_{ij\ell}}.
\label{cosinedihedral formula}
\end{equation}

We use formula (\ref{cosinedihedral formula}) for the cosine of the dihedral
angle and the following expression for the volume of simplex $\left\{
i,j,k,\ell\right\}  $
\[
V_{ijk\ell}=\frac{2A_{ijk}A_{ij\ell}\sin\beta_{ijk\ell}}{3\ell_{ij}}
\]
to compute
\begin{align*}
\frac{\partial\beta_{ijk\ell}}{\partial r_{i}}  &  =\frac{2r_{i}r_{j}r_{k}
^{2}r_{\ell}^{2}}{3P_{ijk}P_{ij\ell}V_{ijk\ell}}\left[  -\frac{1}{r_{k}^{2}
}-\frac{1}{r_{\ell}^{2}}-2\frac{r_{j}}{r_{i}}\left(  \frac{1}{r_{i}\,r_{k}
}+\frac{1}{r_{i}\,r_{\ell}}+\frac{1}{r_{k}\,r_{\ell}}\left(  2+\frac{r_{j}
}{r_{i}}\right)  \right)  \right. \\
&  \quad\quad\,\left.  +\left(  \frac{1}{r_{j}}-\frac{1}{r_{i}}\right)
\left(  \frac{2}{r_{i}}+\frac{1}{r_{k}}+\frac{1}{r_{\ell}}\right)  \right] \\
\frac{\partial\beta_{ijk\ell}}{\partial r_{k}}  &  =\frac{r_{i}^{2}r_{j}
^{2}r_{\ell}}{3P_{ijk}V_{ijk\ell}}\left[  \left(  \frac{1}{r_{i}}+\frac
{1}{r_{j}}\right)  \left(  \frac{1}{r_{i}}+\frac{1}{r_{j}}+\frac{1}{r_{k}
}-\frac{1}{r_{\ell}}\right)  \right]  .
\end{align*}
Now compute the evolution of the solid angles using the formula for the area
of a spherical triangle, $\alpha_{ijk\ell}=\beta_{ijk\ell}+\beta_{ikj\ell
}+\beta_{i\ell jk}-\pi$. We get
\begin{align}
\frac{\partial\alpha_{ijk\ell}}{\partial r_{i}}  &  =-\frac{8r_{j}^{2}
r_{k}^{2}r_{\ell}^{2}}{3P_{ijk}P_{ij\ell}P_{ik\ell}V_{ijk\ell}}\left[  \left(
\frac{2}{r_{i}}+\frac{1}{r_{j}}+\frac{1}{r_{k}}+\frac{1}{r_{\ell}}\right)
\right. \label{partialairi}\\
&  \quad\quad+\frac{r_{j}}{r_{i}}\left(  \frac{1}{r_{i}}+\frac{1}{r_{k}}
+\frac{1}{r_{\ell}}\right)  +\frac{r_{k}}{r_{i}}\left(  \frac{1}{r_{i}}
+\frac{1}{r_{j}}+\frac{1}{r_{\ell}}\right) \nonumber\\
&  \quad\quad+\left.  \frac{r_{\ell}}{r_{i}}\left(  \frac{1}{r_{i}}+\frac
{1}{r_{j}}+\frac{1}{r_{k}}\right)  +\left(  2r_{i}+r_{j}+r_{k}+r_{\ell
}\right)  Q_{ijk\ell}\right] \nonumber
\end{align}
which we see is always negative if the tetrahedron is nondegenerate, i.e.
$Q_{ijk\ell}>0$. The other partial derivatives look like
\begin{align}
\frac{\partial\alpha_{ijk\ell}}{\partial r_{j}}  &  =\frac{4r_{i}r_{j}
r_{k}^{2}r_{\ell}^{2}}{3P_{ijk}P_{ij\ell}V_{ijk\ell}}\left(  \frac{1}{r_{i}
}\left(  \frac{1}{r_{j}}+\frac{1}{r_{k}}+\frac{1}{r_{\ell}}\right)  +\frac
{1}{r_{j}}\left(  \frac{1}{r_{i}}+\frac{1}{r_{k}}+\frac{1}{r_{\ell}}\right)
\right. \label{partialairj}\\
&  \quad\quad\quad\quad-\left.  \left(  \frac{1}{r_{k}}-\frac{1}{r_{\ell}
}\right)  ^{2}\right) \nonumber
\end{align}
which we would like to say is positive, but is not always (although in the
case of most \textquotedblleft good\textquotedblright\ tetrahedra, it is positive).

Finally we sum cyclically in the last three indices and find
\[
\frac{d}{dt}\left(  \beta_{ijk\ell}+\beta_{ikj\ell}+\beta_{i\ell jk}\right)
=\Omega_{ijk\ell}\left(  K_{i}-K_{j}\right)  +\Omega_{ikj\ell}\left(
K_{i}-K_{k}\right)  +\Omega_{i\ell jk}\left(  K_{i}-K_{\ell}\right)
\]
with
\begin{align*}
\Omega_{ijk\ell} &  =\frac{4r_{i}r_{j}^{2}r_{k}^{2}r_{\ell}^{2}}
{3P_{ijk}P_{ij\ell}V_{ijk\ell}}\left(  \frac{1}{r_{i}}\left(  \frac{1}{r_{j}
}+\frac{1}{r_{k}}+\frac{1}{r_{\ell}}\right)  +\frac{1}{r_{j}}\left(  \frac
{1}{r_{i}}+\frac{1}{r_{k}}+\frac{1}{r_{\ell}}\right)  \right.  \\
&  \quad\quad\quad\quad-\left.  \left(  \frac{1}{r_{k}}-\frac{1}{r_{\ell}
}\right)  ^{2}\right)  .
\end{align*}
Thus the evolution of curvature is
\begin{align*}
\frac{dK_{i}}{dt} &  =-\sum_{\left\{  i,j,k,\ell\right\}  \in\mathcal{S}_{3}
}\left(  \frac{d}{dt}\beta_{ijk\ell}+\frac{d}{dt}\beta_{ikj\ell}+\frac{d}
{dt}\beta_{i\ell jk}\right)  .\\
&  =\sum_{\left\{  i,j,k,\ell\right\}  \in\mathcal{S}_{3}}\left[
\Omega_{ijk\ell}\left(  K_{j}-K_{i}\right)  +\Omega_{ikj\ell}\left(
K_{k}-K_{i}\right)  +\Omega_{i\ell jk}\left(  K_{\ell}-K_{i}\right)  \right]
.
\end{align*}
We can define our Laplacian as
\begin{equation}
\left(  \Delta f\right)  _{i}=\sum_{\left\{  i,j,k,\ell\right\}
\in\mathcal{S}_{3}}\left[  \Omega_{ijk\ell}\left(  f_{j}-f_{i}\right)
+\Omega_{ikj\ell}\left(  f_{k}-f_{i}\right)  +\Omega_{i\ell jk}\left(
f_{\ell}-f_{i}\right)  \right]  \label{our laplacian}
\end{equation}
in order to write the evolution of curvature as
\[
\frac{dK_{i}}{dt}=\left(  \triangle K\right)  _{i}.
\]
Unfortunately the coefficients $\Omega_{ijk\ell}$ are not always positive. We
notice that $\Delta$ is self-adjoint with respect to the inner product
\[
\left\langle f,g\right\rangle =\sum_{i\in\mathcal{S}_{0}}f_{i}g_{i}r_{i}
\]
and satisfies
\[
\sum_{i\in\mathcal{S}_{0}}\Delta f_{i}r_{i}=0,
\]
which is analogous in the smooth category to $\int\Delta f\ dx=0$ by the
divergence theorem.

\section{Combinatorial Laplacians}

In this section we investigate the Laplacian defined in (\ref{our laplacian})
and similar operators. Z. He looked at variations of the type we are studying,
that is
\[
\frac{dr_{i}}{dt}=-f_{i}r_{i}
\]
in the two dimensions \cite{he}. He then derived a curvature evolution which
involved a combinatorial Laplacian. The work is similar to the smooth
derivation of the evolution of curvature under a conformal flow on a
Riemannian manifold given by
\[
\frac{\partial}{\partial t}g_{ij}=-fg_{ij}.
\]
This was studied by Hamilton in two dimensions for the Ricci flow, and the
result is
\[
\frac{\partial R}{\partial t}=\Delta_{g}f+Rf
\]
where $\Delta_{g}$ is the Laplacian with respect to the metric $g$ (see
\cite{hamilton:surfaces}).

\begin{figure}
[ptb]
\begin{center}
\includegraphics[
natheight=3.555200in,
natwidth=3.555200in,
height=1.5544in,
width=1.5544in
]
{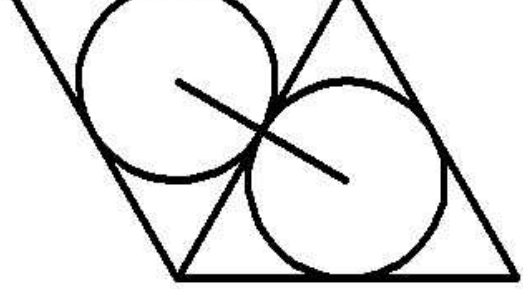}
\caption{Two triangles with inscribed circles and dual edge.}
\label{dueledge}
\end{center}
\end{figure}

Z. He's Laplacian is the same Laplacian derived by Chow and Luo in
\cite{chowluo}, which can be written as
\begin{equation}
\triangle f_{i}=\sum_{\left\{  i,j\right\}  \in\mathcal{S}_{1}}\frac{\ell
_{ij}^{\ast}}{\ell_{ij}}\left(  f_{j}-f_{i}\right)
\label{chow-luo geometric fmla}
\end{equation}
using the relation
\begin{equation}
\frac{\partial\theta_{ijk}}{\partial r_{j}}r_{j}=\frac{r_{ijk}}{\ell_{ij}},
\label{partial derivative is ratio of dual lengths}
\end{equation}
where $\theta_{ijk}$ is the angle at vertex $i$ in triangle $\left\{
i,j,k\right\}  $ and $r_{ijk}$ is the length of the radius of the circle
inscribed in triangle $\left\{  i,j,k\right\}  .$ Hence to get
(\ref{chow-luo geometric fmla}) we simply add
\[
\frac{\partial\theta_{ijk}}{\partial r_{j}}r_{j}+\frac{\partial\theta_{ij\ell
}}{\partial r_{j}}r_{j}=\frac{\ell_{ij}^{\ast}}{\ell_{ij}}
\]
where $\ell_{ij}^{\ast}=r_{ijk}+r_{ij\ell}$ is the length of the dual edge
(see Figure \ref{dueledge}).
The dual vertex $\bigstar\left\{  i,j,k\right\}  $ to the triangle $\left\{
i,j,k\right\}  $ is the center of the inscribed circle, while the edge
$\bigstar\left\{  i,j\right\}  $ dual to $\left\{  i,j\right\}  $ is the edge
which goes from the dual vertex $\bigstar\left\{  i,j,k\right\}  $ to the dual
vertex $\bigstar\left\{  i,j,\ell\right\}  $ which is perpendicular to
$\left\{  i,j\right\}  .$ This Laplacian is similar to the Laplacians found in
the image processing literature which we shall now describe.

Combinatorial Laplacians on piecewise linear surfaces are used quite a bit in
image processing, for instance \cite{meyeretal}, \cite{yauetal}. There is a
very clear description by Hirani \cite{hiranithesis} which defines the
Laplace-Beltrami operator on functions defined at the vertices as
\begin{align}
\triangle f_{i}  &  =\frac{1}{\left\vert \bigstar\left\{  i\right\}
\right\vert }\sum_{\left\{  i,j\right\}  \in\mathcal{S}_{1}}\frac{\left\vert
\bigstar\left\{  i,j\right\}  \right\vert }{\left\vert \left\{  i,j\right\}
\right\vert }\left(  f_{j}-f_{i}\right) \label{laplacian}\\
&  =\frac{1}{V_{i}^{\ast}}\sum_{\left\{  i,j\right\}  \in\mathcal{S}_{1}}
\frac{\ell_{ij}^{\ast}}{\ell_{ij}}\left(  f_{j}-f_{i}\right) \nonumber
\end{align}
where $\bigstar\sigma^{k}$ is the $\left(  n-k\right)  $-dimensional dual of
the $k$-dimensional simplex $\sigma^{k}$ and $\left\vert \sigma^{k}\right\vert
$ is the $k$-dimensional volume of the simplex. The second line uses our
notation, where the dual of an edge $\ell_{ij}^{\ast}$ and the volume of a
vertex $V_{i}^{\ast}$ are defined appropriately. Note that Z. He's Laplacian
is exactly this, except for the volume factor in front, where duality comes
from the inscribed circles. We also note that in most of the image processing
literature, the two-dimensional dual comes from the center of the
circumscribed circle instead of the inscribed circle described here.

It is interesting to note the geometric justification for the formula
(\ref{laplacian}). If we consider the integral of the Laplace-Beltrami
operator and use Stokes' theorem, we find
\[
\int_{U}\triangle f~dV=\int_{\partial U}\frac{df}{dn}dS
\]
where $\frac{df}{dn}$ is the normal derivative and $dS$ is the surface
measure. We easily see that if we take $U$ to be $\bigstar\left\{  i\right\}
$ then the normal derivative is
\[
\frac{f_{j}-f_{i}}{\ell_{ij}},
\]
$\ell_{ij}^{\ast}$ is the surface measure, and $V_{i}^{\ast}$ is the volume measure.

\begin{figure}
[ptb]
\begin{center}
\includegraphics[
natheight=2.701900in,
natwidth=5.217500in,
height=1.7868in,
width=3.4334in
]
{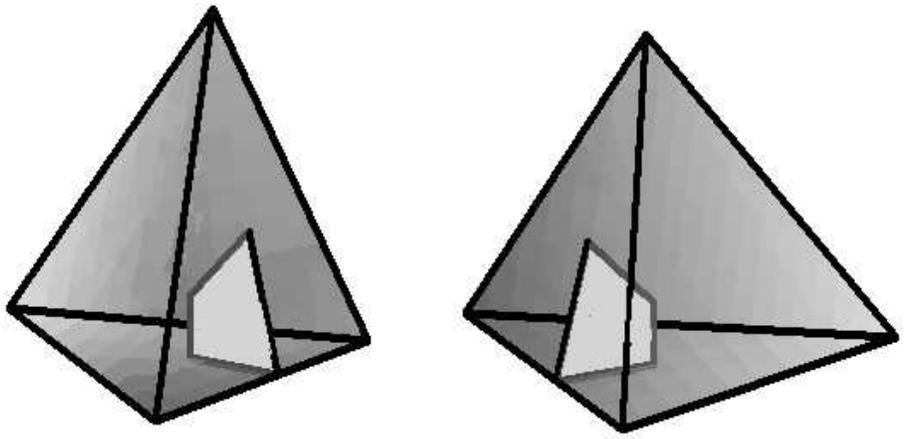}
\caption{Two views of the dual to an edge.}
\label{tetraduals}
\end{center}
\end{figure}

Recall the definition (\ref{our laplacian}) of the Laplacian we gave in the
previous section. This Laplacian would be related to the Laplacian in this
section (\ref{laplacian}) if we had an analogue of
(\ref{partial derivative is ratio of dual lengths}). The dual $\bigstar
\left\{  i,j\right\}  $ to an edge $\left\{  i,j\right\}  $ is a surface which
goes through $\bigstar T$ for any tetrahedron $T$ containing $\left\{
i,j\right\}  $ and is perpendicular to $\left\{  i,j\right\}  .$ Two pictures
of the piece of $\bigstar\left\{  i,j\right\}  $ in one tetrahedron can be
seen in Figure \ref{tetraduals}.
We shall denote the area shown in one tetrahedron by $A_{ijk\ell},$ where the
region is perpendicular to the edge $\left\{  i,j\right\}  .$ The region whose
area is $A_{ijk\ell}$ has four sides:

\begin{itemize}
\item the radius $r_{ijk}$ of the circle inscribed in the face $\left\{
i,j,k\right\}  $ which intersects $\left\{  i,j\right\}  ,$

\item the radius $r_{ij\ell}$ of the circle inscribed in the face $\left\{
i,j,\ell\right\}  $ which intersects $\left\{  i,j\right\}  ,$

\item the line from the center of the circumscripted sphere to the face
$\left\{  i,j,k\right\}  $ which is perpendicular to the plane determined by
$\left\{  i,j,k\right\}  ,$ and

\item the line from the center of the circumscripted sphere to the face
$\left\{  i,j,\ell\right\}  $ which is perpendicular to the plane determined
by $\left\{  i,j,\ell\right\}  .$
\end{itemize}

\noindent It can also be decomposed as two right triangles, each having
hypotenuse coming from the radius $r_{ijk\ell}$ of the circumscripted sphere
which intersects the edge $\left\{  i,j\right\}  .$ One triangle has a leg
coming from the radius $r_{ijk}$ of the circle inscribed in $\left\{
i,j,k\right\}  $ and the other has a leg coming from the radius $r_{ij\ell}$
of the circle inscribed in $\left\{  i,j,\ell\right\}  $. The legs and
hypotenuse meet at the same point in $\left\{  i,j\right\}  .$ We shall call
the lengths of the other respective legs $h_{ijk,\ell}$ and $h_{ij\ell,k}$ to
denote that they are heights of tetrahedra which make up the tetrahedron
$\left\{  i,j,k,\ell\right\}  ;$ that is, if we let $c$ denote the center of
the circumscripted sphere, $h_{ijk,\ell}$ is the height of tetrahedron
$\left\{  i,j,k,c\right\}  $ with base $\left\{  i,j,k\right\}  .$ Thus the
volume of tetrahedron $\left\{  i,j,k,\ell\right\}  $ is decomposed as
\begin{equation}
V_{ijk\ell}=\sum_{\left\{  i,j,k\right\}  \in\mathcal{T}_{2}}\frac{1}
{3}A_{ijk}h_{ijk,\ell}. \label{volume decomposition of tetrahedra}
\end{equation}
If the center of the circumscripted sphere is outside the tetrahedron, we
define $h_{ijk,\ell}$ to be negative if the center is on the opposite side of
the plane determined by $\left\{  i,j,k\right\}  $ from the tetrahedron. This
way, (\ref{volume decomposition of tetrahedra}) is still satisfied. The
$h_{ijk,\ell}$ are symmetric in the first three indices, and the fourth
indicates to which tetrahedron the $h$ corresponds. The area $A_{ijk\ell}$ is
then computed as
\[
A_{ijk\ell}=\frac{1}{2}h_{ijk,\ell}r_{ijk}+\frac{1}{2}h_{ij\ell,k}r_{ij\ell},
\]
where this may be negative.

We can now show the analogue of
(\ref{partial derivative is ratio of dual lengths}).

\begin{lem}
\[
\frac{\partial\alpha_{ijk\ell}}{\partial r_{j}}r_{i}r_{j}=\frac{A_{ijk\ell}
}{\ell_{ij}}
\]
where $A_{ijk\ell}$ is the (signed) area of the dual region to the side
$\left\{  i,j\right\}  $ in the tetrahedron $\left\{  i,j,k,\ell\right\}  $ as
described above.
\end{lem}

\begin{pf}
We simply compute. Recall that $P_{ijk}$ is the perimeter of $\left\{
i,j,k\right\}  ,$ $A_{ijk}$ is the area of $\left\{  i,j,k\right\}  ,$
$r_{ijk}$ is the radius of the circle inscribed in $\left\{  i,j,k\right\}  ,$
$r_{ijk\ell}$ is the radius of the sphere circumscripted by $\left\{
i,j,k,\ell\right\}  ,$ and $h_{ijk,\ell}$ is signed height of the tetrahedron
defined by the center of the circumscripted sphere in $\left\{  i,j,k,\ell
\right\}  $ and $i,$ $j,$ and $k,$ with base $\left\{  i,j,k\right\}  .$ The
sign of $h_{ijk,\ell}$ is defined so that
\[
3V_{ijk\ell}=h_{ijk,\ell}A_{ijk}+h_{ij\ell,k}A_{ij\ell}+h_{ik\ell,j}A_{ik\ell
}+h_{jk\ell,i}A_{jk\ell}.
\]
We have the following relations:
\begin{align*}
A_{ijk}  &  =\frac{1}{2}r_{ijk}P_{ijk}\\
r_{ijk}  &  =\frac{2A_{ijk}}{P_{ijk}}=\sqrt{\frac{r_{i}r_{j}r_{k}}{r_{i}
+r_{j}+r_{k}}}.
\end{align*}
We also know that
\[
\frac{4}{r_{ijk\ell}^{2}}=Q_{ijk\ell}=\left(  \frac{1}{r_{i}}+\frac{1}{r_{j}
}+\frac{1}{r_{k}}+\frac{1}{r_{\ell}}\right)  ^{2}-2\left(  \frac{1}{r_{i}^{2}
}+\frac{1}{r_{j}^{2}}+\frac{1}{r_{k}^{2}}+\frac{1}{r_{\ell}^{2}}\right)  .
\]
The square of the height $h_{ijk,\ell}$ can be computed using the Pythagorean
theorem:
\begin{align*}
h_{ijk,\ell}^{2}  &  =r_{ijk\ell}^{2}-r_{ijk}^{2}\\
&  =\frac{r_{ijk\ell}^{2}r_{ijk}^{2}}{4}\left(  4\frac{r_{i}+r_{j}+r_{k}
}{r_{i}r_{j}r_{k}}-\left(  \left(  \frac{1}{r_{i}}+\frac{1}{r_{j}}+\frac
{1}{r_{k}}+\frac{1}{r_{\ell}}\right)  ^{2}-2\left(  \frac{1}{r_{i}^{2}}
+\frac{1}{r_{j}^{2}}+\frac{1}{r_{k}^{2}}+\frac{1}{r_{\ell}^{2}}\right)
\right)  \right) \\
&  =\frac{r_{ijk\ell}^{2}r_{ijk}^{2}}{4}\left(  \frac{1}{r_{i}}+\frac{1}
{r_{j}}+\frac{1}{r_{k}}-\frac{1}{r_{\ell}}\right)  ^{2}.
\end{align*}
In order to determine $h_{ijk,\ell}$ correctly, we notice that if the center
of the circumscripted sphere is on face $\left\{  i,j,k\right\}  $ then
\[
r_{ijk}^{2}=r_{ijk\ell}^{2}
\]
and hence
\begin{align*}
4\frac{r_{i}+r_{j}+r_{k}}{r_{i}r_{j}r_{k}}  &  =\left(  \frac{1}{r_{i}}
+\frac{1}{r_{j}}+\frac{1}{r_{k}}+\frac{1}{r_{\ell}}\right)  ^{2}-2\left(
\frac{1}{r_{i}^{2}}+\frac{1}{r_{j}^{2}}+\frac{1}{r_{k}^{2}}+\frac{1}{r_{\ell
}^{2}}\right) \\
0  &  =\left(  \frac{1}{r_{i}}+\frac{1}{r_{j}}+\frac{1}{r_{k}}-\frac
{1}{r_{\ell}}\right)  ^{2}.
\end{align*}
We easily see that if the vector from the center to the side $\left\{
i,j,k\right\}  $ is in the same direction as the outward pointing normal of
side $\left\{  i,j,k\right\}  $ then
\[
\frac{1}{r_{i}}+\frac{1}{r_{j}}+\frac{1}{r_{k}}>\frac{1}{r_{\ell}}
\]
and if the vector is in the opposite direction than the outward pointing
normal then
\[
\frac{1}{r_{i}}+\frac{1}{r_{j}}+\frac{1}{r_{k}}<\frac{1}{r_{\ell}}.
\]
Hence the signed height $h_{ijk,\ell}$ is
\[
h_{ijk,\ell}=\frac{r_{ijk\ell}r_{ijk}}{2}\left(  \frac{1}{r_{i}}+\frac
{1}{r_{j}}+\frac{1}{r_{k}}-\frac{1}{r_{\ell}}\right)  ,
\]
which may be negative. Notice that this assures that
\[
3V_{ijk\ell}=h_{ijk,\ell}A_{ijk}+h_{ij\ell,k}A_{ij\ell}+h_{ik\ell,j}A_{ik\ell
}+h_{jk\ell,i}A_{jk\ell},
\]
where $V_{ijk\ell}$ is the volume of $\left\{  i,j,k,\ell\right\}  .$ The dual
area is computed to be
\begin{align*}
A_{ijk\ell}  &  =\frac{1}{2}h_{ijk,\ell}r_{ijk}+\frac{1}{2}h_{ij\ell
,k}r_{ij\ell}\\
&  =\frac{r_{ijk\ell}}{4}\left(  \frac{4A_{ijk}^{2}}{P_{ijk}^{2}}\left(
\frac{1}{r_{i}}+\frac{1}{r_{j}}+\frac{1}{r_{k}}-\frac{1}{r_{\ell}}\right)
+\frac{4A_{ij\ell}^{2}}{P_{ij\ell}^{2}}\left(  \frac{1}{r_{i}}+\frac{1}{r_{j}
}-\frac{1}{r_{k}}+\frac{1}{r_{\ell}}\right)  \right) \\
&  =\frac{4r_{i}^{2}r_{j}^{2}r_{k}^{2}r_{\ell}^{2}\ell_{ij}}{3V_{ijk\ell
}P_{ijk}P_{ij\ell}}\left(  \frac{1}{r_{i}}\left(  \frac{1}{r_{j}}+\frac
{1}{r_{k}}+\frac{1}{r_{\ell}}\right)  +\frac{1}{r_{j}}\left(  \frac{1}{r_{i}
}+\frac{1}{r_{k}}+\frac{1}{r_{\ell}}\right)  -\left(  \frac{1}{r_{k}}-\frac
{1}{r_{\ell}}\right)  ^{2}\right)  .
\end{align*}
\qed

\end{pf}

\begin{cor}
We have that
\[
\triangle f_{i}=\frac{1}{r_{i}}\sum_{\left\{  i,j\right\}  \in\mathcal{S}_{1}
}\frac{\ell_{ij}^{\ast}}{\ell_{ij}}\left(  f_{j}-f_{i}\right)  ,
\]
where $\ell_{ij}^{\ast}$ is the area dual to the side $\left\{  i,j\right\}
.$
\end{cor}

\begin{pf}
This follows from the fact that the dual area $\ell_{ij}^{\ast}$ is simply
\[
\ell_{ij}^{\ast}=\sum_{\left\{  i,j,k,\ell\right\}  \in\mathcal{S}_{3}
}A_{ijk\ell},
\]
where the sum is over $k$ and $\ell,$ that is, all tetrahedra incident on the
edge $\left\{  i,j\right\}  .$ \qed
\end{pf}

Note the similarity to Hirani's definition (\ref{laplacian}). Also note that
since $\ell_{ij}^{\ast}$ may be negative, this is not always a Laplacian on
graphs in the usual sense (see, for instance, \cite{chung}). Still, we can
prove the maximum principle in more general circumstances, as seen in
\cite{glickensteinmaxprinciple}.

\section{Convergence to constant curvature}

In this section we shall show that if the solution exists for all time, the
curvatures converge to a constant. We restrict ourselves to well behaving solutions.

\begin{defn}
\label{nonsingular def}$\left\{  r_{i}\left(  t\right)  \right\}
_{i\in\mathcal{S}_{0}}$ is a \emph{nonsingular solution }of the combinatorial
Yamabe flow if there exists $\delta>0$ such that for each for $t\in
\lbrack0,\infty)$ it satisfies (\ref{combinatorial yamabe flow}), $Q_{ijk\ell
}>0$ for all $\left\{  i,j,k,\ell\right\}  \in\mathcal{S}_{3},$ and
\[
\frac{r_{i}}{\sum_{j\in\mathcal{S}_{0}}r_{j}}\geq\delta
\]
for each $i\in\mathcal{S}_{0}.$
\end{defn}

Consider the average scalar curvature
\[
k\doteqdot\frac{\sum_{i\in\mathcal{S}_{0}}K_{i}r_{i}}{\sum_{i\in
\mathcal{S}_{0}}r_{i}}.
\]
This can be thought of as an analogue of the average scalar curvature
function
\[
\frac{\int_{M}R~dV}{\int_{M}dV}
\]
on a Riemannian manifold. Note that this curvature really is an average in the
sense that
\[
K_{\min}\leq k\leq K_{\max}
\]
if $K_{\min}$ and $K_{\max}$ are the minimal and maximal curvatures.

\begin{prop}
\label{curvature converges to constant}For any nonsingular solution $\left\{
r_{i}\left(  t\right)  \right\}  _{i\in\mathcal{S}_{0}}$ there is a number
$k\left(  \infty\right)  $ such that the average scalar curvature $k\left(
t\right)  $ and all of the curvatures $K_{i}\left(  t\right)  $ converge to
$k\left(  \infty\right)  $ as $t\rightarrow\infty.$
\end{prop}

\begin{pf}
It is easy to see that $k$ is decreasing along the flow by a direct
computation:
\begin{align*}
\frac{dk}{dt}  &  =-\frac{\sum_{i\in\mathcal{S}_{0}}K_{i}^{2}r_{i}}{\sum
_{i\in\mathcal{S}_{0}}r_{i}}+\frac{\left(  \sum_{i\in\mathcal{S}_{0}}
K_{i}r_{i}\right)  ^{2}}{\left(  \sum_{i\in\mathcal{S}_{0}}r_{i}\right)  ^{2}
}\\
&  =-\frac{\sum_{i\in\mathcal{S}_{0}}\sum_{j\in\mathcal{S}_{0}}\left(
K_{i}^{2}r_{i}r_{j}-K_{i}K_{j}r_{i}r_{j}\right)  }{\left(  \sum_{i\in
\mathcal{S}_{0}}r_{i}\right)  ^{2}}.
\end{align*}
Rearranging terms we get
\begin{equation}
\frac{dk}{dt}=-\frac{\sum_{i\in\mathcal{S}_{0}}\sum_{j\in\mathcal{S}_{0}
}\left(  K_{i}-K_{j}\right)  ^{2}r_{i}r_{j}}{\left(  \sum_{i\in\mathcal{S}
_{0}}r_{i}\right)  ^{2}}. \label{deriv of average curvature}
\end{equation}
Furthermore,
\[
K_{\min}\leq k
\]
where $K_{\min}$ is the minimum of the curvature, which is bounded below by
$4\pi-2\pi d_{\max}$ if $d_{\max}$ is the maximum number of tetrahedra
incident on any one vertex. Thus $k$ is decreasing and bounded below, so it
must converge to a limit $k\left(  \infty\right)  .$ Moreover, the time
derivative of $k$ must go to zero or $k$ would not be bounded below. Hence
$\frac{dk}{dt}$ converges to zero. By formula
(\ref{deriv of average curvature}) and the fact that
\[
\frac{r_{i}}{\sum r_{j}}\geq\delta
\]
we see that $\left(  K_{i}-K_{j}\right)  ^{2}\rightarrow0$ for all pairs of
vertices. Hence the curvatures becomes constant; this constant must be
$k\left(  \infty\right)  .$ \qed

\end{pf}

\section{Long term existence}

In this section we will classify the possible long term behavior.

\begin{prop}
All solutions to the combinatorial Yamabe flow on a maximal time interval
$[0,T)$ must fit into one of the following categories:

\begin{itemize}
\item It is nonsingular (see Definition \ref{nonsingular def}).

\item $T=\infty$ and for some $i\in\mathcal{S}_{0},$
\[
\frac{r_{i}}{\sum_{j\in\mathcal{S}_{0}}r_{j}}\rightarrow0
\]
as $t\rightarrow\infty.$

\item $T<\infty$ and for some $\left\{  i,j,k,\ell\right\}  \in\mathcal{S}
_{3},$ $Q_{ijk\ell}\rightarrow0$ as $t\nearrow T.$
\end{itemize}
\end{prop}

Solutions with an infinite time interval are covered, and so we need only look
at finite time singularities. \textit{A priori}, the following may also happen:

\begin{itemize}
\item there exists $i\in\mathcal{S}_{0}$ such that $r_{i}\left(  t\right)
\rightarrow0$ as $t\nearrow T$

\item there exists $i\in\mathcal{S}_{0}$ such that $r_{i}\left(  t\right)
\rightarrow\infty$ as $t\nearrow T.$
\end{itemize}

We can consider $L_{i}\left(  t\right)  =\log r_{i}\left(  t\right)  .$ Notice
that
\[
\frac{dL_{i}}{dt}=-K_{i}.
\]
If $L_{i}\rightarrow\pm\infty$ in finite time, then $\frac{dL_{i}}{dt}
=K_{i}\rightarrow\pm\infty$ in finite time. This is impossible, though, since
\[
4\pi-2\pi d_{\max}\leq K_{i}\leq4\pi.
\]
So finite time singularities occur only because $Q_{ijk\ell}\rightarrow0.$ The
case of $T=\infty$ and there exists $i$ such that
\[
\frac{r_{i}}{\sum r_{j}}\rightarrow0
\]
is somehow analogous to collapse in the smooth case. Hence Proposition
\ref{curvature converges to constant} is an analogue of Hamilton's theorem on
nonsingular solutions \cite{hamilton:nonsingular}.

\section{Further Remarks}

It is quite easy to implement the combinatorial Yamabe flow numerically. With
the help of F.H. Lutz's work (see \cite{lutzthesis} and
\cite{lutzmanifoldpage}) on small triangulations of manifolds, the author has
been able to run examples of the combinatorial Yamabe flow on manifolds
homeomorphic to the sphere, torus, $S^{2}\times S^{1},$ and $S^{2}
\tilde{\times}S^{1}.$ In each case the combinatorial Yamabe flow found a
constant curvature metric. However, these were all done with relatively small
triangulations so many of the possible degeneracies coming from a poor
triangulation cannot occur.

As noted in Remark \ref{vertex transitive remark}, any metric which is vertex
transitive, that is, has the same number of tetrahedra incident on each
vertex, will have a constant curvature metric where all $r_{i}$ are equal. It
is not clear that this is the only constant curvature metric in the conformal
class. Numerical simulation has been unable to produce any examples of other
constant curvature metrics, however. In the case of the torus, the metric
where all $r_{i}$ are equal is usually not one where the curvatures are zero.
In all simulations of small triangulations of the torus tried, the
combinatorial Yamabe flow converges to a positively curved constant curvature metric.

Constant curvature as indicated in this paper is a weak condition. A metric
such that the sum of the dihedral angles along each edge is equal to $2\pi$ is
a Euclidean structure on a manifold, but a metric where the solid angles equal
$4\pi$ is not necessarily. The difference is something like the difference
between constant sectional curvatures of zero and constant scalar curvature of
zero; it is possible to have a metric where the latter is true but not the former.

\section*{Acknowledgements}

The author would like to thank Ben Chow for introducing the combinatorial
Yamabe flow and for all his help. The author would also like to thank Feng Luo
and Mauro Carfora for useful conversations.

\end{document}